\documentclass[11pt]{article}
\usepackage{amsmath}
\usepackage{amssymb}
\usepackage{amsthm}

\usepackage{url}
\newcommand\nonu{\nonumber}
\newcommand\dstyle\displaystyle
\newcommand\sLP{\\[\smallskipamount]}
\newcommand\mPP{\\[\medskipamount]\indent}
\newcommand\bPP{\\[\bigskipamount]\indent}
\newcommand\al\alpha
\newcommand\be\beta
\newcommand\ga\gamma
\newcommand\de\delta
\newcommand\tha\theta
\newcommand\la\lambda
\newcommand\si\sigma
\newcommand\half{\frac12}
\newcommand\thalf{\tfrac12}
\newcommand\iy\infty
\newcommand\LHS{left-hand side}
\newcommand\RHS{right-hand side}
\newcommand{\hyp}[5]{\,\mbox{}_{#1}F_{#2}\!\left(
  \genfrac{}{}{0pt}{}{#3}{#4};#5\right)}
\newcommand{\qhyp}[5]{\,\mbox{}_{#1}\phi_{#2}\!\left(
  \genfrac{}{}{0pt}{}{#3}{#4};#5\right)}
\newcommand\qbinom[3]{\genfrac[]{0pt}0{#1}{#2}_{#3}}
\newcommand{\dup}{\textup{d}}
\newcommand{\eup}{\textup{e}}
\newcommand{\iup}{\mkern1mu\textup{i}\mkern1mu}

\numberwithin{equation}{section}
\newtheorem{theorem}{Theorem}[section]
\newcommand\Proof{\noindent{\bf Proof}\quad}
\newtheorem{Remark}[theorem]{Remark}
\newenvironment{remark}{\begin{Remark}\rm}{\end{Remark}}

\begin{document}
\title{Dual addition formulas: the case of continuous\\
$q$-ultraspherical and $q$-Hermite polynomials}
\author{Tom H. Koornwinder}
\date{This paper is dedicated to the memory of Dick Askey}
\maketitle
\begin{abstract}
We settle the dual addition formula for continuous
$q$-ultraspherical polynomials as an expansion in terms of
special $q$-Racah polynomials for which the constant term is
given by the linearization formula for the
continuous $q$-ultraspherical polynomials. In a second proof
we derive the dual addition formula from the Rahman--Verma addition
formula for these polynomials by using the self-duality of the
polynomials.
We also consider the limit case of continuous $q$-Hermite polynomials.
\end{abstract}
\section{Introduction}
In this paper, as a natural continuation of our recent derivation \cite{8} of
the dual addition formula for ultraspherical polynomials, we
derive the dual addition formula for continuous
$q$-ultraspherical polynomials.
We give two different proofs. The first proof is a perfect $q$-analogue
of the derivation in \cite{8}. Every step of the proof yields in
the limit for $q\to1$ the corresponding step in \cite{8}.
The second proof exploits the self-duality of the
continuous $q$-ultraspherical polynomials. Then the dual addition
formula easily follows from the known addition formula \cite{10}
for these polynomials.

Addition formulas are closely related to product formulas.
For instance, the \emph{addition formula for Legendre polynomials}
\cite[(18.18.9)]{11}
\begin{multline}
P_n(\cos\tha_1\cos\tha_2+\sin\tha_1\sin\tha_2\cos\phi)
=P_n(\cos\tha_1)P_n(\cos\tha_2)\\
+2\sum_{k=1}^n \frac{(n-k)!\,(n+k)!}{2^{2k}(n!)^2}\,
(\sin\tha_1)^k P_{n-k}^{(k,k)}(\cos\tha_1)
\,(\sin\tha_2)^k P_{n-k}^{(k,k)}(\cos\tha_2)\cos(k\phi)
\label{70}
\end{multline}
gives the Fourier-cosine expansion of the \LHS\ as a function
of $\phi$. Integration with respect to $\phi$ over $[0,\pi]$
gives the constant term in this expansion, which is
the \emph{product formula for Legendre polynomials} \cite[(18.17.6)]{11}
\begin{equation}
P_n(\cos\tha_1)P_n(\cos\tha_2)=\frac1\pi \int_0^\pi
P_n(\cos\tha_1\cos\tha_2+\sin\tha_1\sin\tha_2\cos\phi)\,\dup\phi.
\label{68}
\end{equation}

Two formulas involving Legendre polynomials $P_n(x)$ (or more generally
some orthogonal polynomials $p_n(x)$) are called \emph{dual} to
each other if the roles of $n$ and $x$ in the second formula are
interchanged in comparison with the first formula. The formula
dual to the product formula \eqref{68} is the
\emph{linearization formula}, which expands the product
$P_\ell(x)P_m(x)$ in terms of Legendre polynomials $P_k(x)$.
This expansion is a sum running from $k=|\ell-m|$ to $k=\ell+m$,
where only terms with $\ell+m-k$ even will occur since
$P_n(-x)=(-1)^n P_n(x)$. The linearization formula for Legendre polynomials is explicitly known
(see \cite[(18.18.22)]{11} for $\la=\half$
together with \cite[(18.7.9)]{11}):
\begin{equation}
P_\ell(x)P_m(x)=\sum_{j=0}^{\min(\ell,m)}
\frac{(\half)_j(\half)_{\ell-j}(\half)_{m-j} (\ell+m-j)!}
{j!\,(\ell-j)!\,(m-j)!\,
(\tfrac32)_{\ell+m-j}}\,\big(2(\ell+m-2j)+1\big)\,P_{\ell+m-2j}(x),
\label{69}
\end{equation}
where $(a)_k$ is the shifted factorial, see below.
Dick Askey, in his lectures at conferences, often
raised the problem to find an addition type formula associated with
\eqref{69} in a similar way as the addition formula \eqref{70} is
associated with the product formula \eqref{68}.
The author finally solved this in \cite{8} by recognizing the coefficient
of $P_{\ell+m-2j}(x)$ in \eqref{69} as the weight of a special
Racah polynomial \cite[(9.2.1)]{2} depending on~$j$, and then finding the
expansion of $P_{\ell+m-2j}(x)$ in terms of these Racah polynomials.
More generally, the same idea worked out well in \cite{8} for
ultraspherical polynonmials.

While \eqref{70}, \eqref{68}, \eqref{69}, and their
generalizations to ultraspherical polynomials, are formulas
established long ago and staying within the realm of classical
orthogonal polynomials, it is remarkable that the dual addition formula
steps out from this and needs Racah polynomials, which live
high up in the Askey scheme.
Parallel to the Askey scheme there is the much larger
$q$-Askey scheme\footnote{See
\url{http://homepage.tudelft.nl/11r49/book.html} for charts of these schemes.}.
Families of orthogonal polynomials in the Askey scheme are limit cases
of families in the $q$-Askey scheme. The continuous
$q$-ultraspherical polynomials form the family which is the $q$-analogue
of the ultraspherical polynomials. Moreover, the $q$-analogues
of \eqref{70}, \eqref{68} and \eqref{69} for these polynomials are
available in the literature. The continuous $q$-ultraspherical
polynomials also have the property of \emph{self-duality}, which is
lost in the limit to $q=1$. This notion means that, for a suitable
function $\si$, an orthogonal polynomial $p_n(x)$ has the property
that $p_n(\si(m))=p_m(\si(n))$ ($m,n=0,1,\ldots$). With all this
material available there is a clear road
to the derivation of the dual addition formula for these polynomials.
\mPP
The contents of the paper are as follows.
Section \ref{94} summarizes the results from \cite{8}
about the dual addition formula for ultraspherical polynomials.
The necessary preliminaries about special orthogonal polynomials in the
$q$-case are given in Section \ref{95}.
The new results of the paper appear in Section \ref{75}.
It contains the two proofs of the dual addition formula for
continuous $q$-ultraspherical polynomials. Finally the limit
case for continuous $q$-Hermite polynomials is considered in Section \ref{62}.
\paragraph{Note}
For definition and notation of ($q$-)shifted factorials and
($q$-)hypergeometric series we follow \cite[\S1.2]{3}.
We will only need terminating series:
\begin{align*}
\hyp rs{-n,a_2,\ldots,a_r}{b_1,\ldots,b_s}z&:=
\sum_{k=0}^n \frac{(-n)_k}{k!}\,\frac{(a_2,\ldots,a_r)_k}
{(b_1,\ldots,b_s)_k}\,z^k,\\
\qhyp rs{q^{-n},a_2,\ldots,a_r}{b_1,\ldots,b_s}{q,z}&:=
\sum_{k=0}^n \frac{(q^{-n};q)_k}{(q;q)_k}\,
\frac{(a_2,\ldots,a_r;q)_k}{(b_1,\ldots,b_s;q)_k}\,
\big((-1)^k q^{\half k(k-1)}\big)^{s-r+1}z^k.
\end{align*}
Here $(b_1,\ldots,b_s)_k:=(b_1)_k\ldots(b_s)_k$ with
$(b)_k:=b(b+1)\cdots(b+k-1)$ the shifted
factorial, and
$(b_1,\ldots,b_s;q)_k:=(b_1;q)_k\ldots(b_s;q)_k$ with
$(b;q)_k:=(1-b)(1-qb)\ldots(1-q^{k-1}b)$ the
$q$-shifted factorial.

For formulas on orthogonal polynomials in the ($q$-)Askey scheme we
will often refer to Chapters 9 and 14 in \cite{2}.
\section{The dual addition formula for ultraspherical polynomials}
\label{94}
Here we summarize the results of \cite{8}. We write \emph{ultraspherical
polynomials} as
\begin{equation}
R_n^\al(x):=\frac{P_n^{(\al,\al)}(x)}{P_n^{(\al,\al)}(1)}=
\frac{C_n^{(\al+\half)}(x)}{C_n^{(\al+\half)}(1)}=
\hyp21{-n,n+2\al+1}{\al+1}{\thalf(1-x)},
\label{90}
\end{equation}
where $C_n^{(\la)}(x)$ is the standard notation \cite[\S9.8.1]{2} for
ultraspherical polynomials and $P_n^{(\al,\be)}(x)$ is a Jacobi polynomial
\cite[\S9.8]{2}.

We will consider \emph{Racah polynomials}
\cite[\S9.2]{2}
\begin{equation}
R_n\big(x(x+\ga+\de+1);\al,\be,\ga,\de\big):=
\hyp43{-n,n+\al+\be+1,-x,x+\ga+\de+1}{\al+1,\be+\de+1,\ga+1}1
\label{91}
\end{equation}
for $\ga=-N-1$, where $N\in\{1,2,\ldots\}$, and for
$n\in\{0,1,\ldots,N\}$. These are orthogonal polynomials on
the finite quadratic set $\{x(x+\ga+\de+1)\mid x\in\{0,1,\ldots,N\}\}$:
\begin{equation*}
\sum_{x=0}^N (R_mR_n)\big(x(x+\ga+\de+1);\al,\be,\ga,\de\big)\,
w_{\al,\be,\ga,\de}(x)=h_{n;\al,\be,\ga,\de}\,\de_{m,n}\quad
(m,n\in\{0,1,\ldots,N\})
\end{equation*}
with
\begin{equation}
w_{\al,\be,\ga,\de}(x)=
\frac{(\al+1)_x (\be+\de+1)_x (\ga+1)_x (\ga+\de+1)_x}
{(-\al+\ga+\de+1)_x (-\be+\ga+1)_x (\de+1)_x\,x!}\,
\frac{\ga+\de+1+2x}{\ga+\de+1},
\label{92}
\end{equation}
\begin{equation*}
\frac{h_{n;\,\al,\be,\ga,\de}}{h_{0;\,\al,\be,\ga,\de}}=
\frac{\al+\be+1}{\al+\be+2n+1}\,
\frac{(\be+1)_n(\al+\be-\ga+1)_n(\al-\de+1)_n\,n!}
{(\al+1)_n(\al+\be+1)_n(\be+\de+1)_n(\ga+1)_n},
\end{equation*}
\begin{equation}
h_{0;\,\al,\be,\ga,\de}=\sum_{x=0}^N w_{\al,\be,\ga,\de}(x)=
\frac{(\al+\be+2)_N(-\de)_N}{(\al-\de+1)_N(\be+1)_N}\qquad(\ga=-N-1).
\label{93}
\end{equation}

The linearization formula for ultraspherical polynomials,
see \cite[(5.7)]{7}, can be written as
\begin{multline}
R_\ell^\al(x) R_m^\al(x)=
\frac{\ell!\,m!}{(2\al+1)_\ell(2\al+1)_m}
\sum_{j=0}^{\min(\ell,m)}\frac{\ell+m+\al+\half-2j}{\al+\half}\\
\times
\frac{(\al+\half)_j(\al+\half)_{\ell-j}(\al+\half)_{m-j}
(2\al+1)_{\ell+m-j}}
{j!\,(\ell-j)!\,(m-j)!\,(\al+\frac32)_{\ell+m-j}}\,
R_{\ell+m-2j}^\al(x).
\label{22}
\end{multline}
Assume that $\al>-\half$ and that, without loss of generality,
$\ell\ge m$.
By \eqref{92} and \eqref{93}
formula \eqref{22} can be rewritten as
\begin{equation}
R_\ell^\al(x) R_m^\al(x)=
\sum_{j=0}^m
\frac{w_{\al-\half,\al-\half,-m-1,-\ell-\al-\half}(j)}
{h_{0;\,\al-\half,\al-\half,-m-1,-\ell-\al-\half}}\,
R_{\ell+m-2j}^\al(x)\quad(\ell\ge m).
\label{27}
\end{equation}
This can be considered as giving the constant term
of an expansion of $R_{\ell+m-2j}^{(\al,\al)}(x)$ as a function of $j$
in terms of the following special case of Racah polynomials \eqref{91}:
\begin{multline*}
R_n\big(j(j-\ell-m-\al-\thalf);\,
\al-\thalf,\al-\thalf,-m-1,-\ell-\al-\thalf\big)\\
=\hyp43{-n,n+2\al,-j,j-\ell-m-\al-\thalf}{\al+\thalf,-\ell,-m}1.
\end{multline*}
The full expansion
is the \emph{dual addition formula for ultraspherical polynomials}\,:
\begin{multline}
R_{\ell+m-2j}^\al(x)
=\sum_{k=0}^{\min(l,m)}\frac{\al+k}{\al+\thalf k}\,
\frac{(-\ell)_k (-m)_k (2\al+1)_k}{2^{2k}(\al+1)_k^2\,k!}\,
(x^2-1)^k\,
R_{\ell-k}^{\al+k}(x)\,R_{m-k}^{\al+k}(x)\\
\times R_k\big(j(j-\ell-m-\al-\thalf);\,
\al-\thalf,\al-\thalf,-m-1,-\ell-\al-\thalf\big),
\quad j\in\{0,1,\ldots,m\}.
\label{29}
\end{multline}
For $j=0$ this becomes
\begin{equation}
R_{\ell+m}^\al(x)
=\sum_{k=0}^{\min(l,m)}\frac{\al+k}{\al+\thalf k}\,
\frac{(-\ell)_k (-m)_k (2\al+1)_k}{2^{2k}(\al+1)_k^2\,k!}\,
(x^2-1)^k\,
R_{\ell-k}^{\al+k}(x)\,R_{m-k}^{\al+k}(x).
\label{97}
\end{equation}
Formula \eqref{97} was first given by Carlitz \cite[(3)]{45}.
It can be rewritten as a matrix decomposition $S=LDU$ with $S$
symmetric, $L$ lower triangular, its transpose $U=L^\intercal$ upper
triangular and $D$ diagonal. 
Cagliero \& Koornwinder \cite[Theorem 4.1 for $\al=\be$]{38} earlier gave
the inverse of the matrix $L$.
\section{Some $q$-hypergeometric orthogonal polynomials}
\label{95}
\subsection{Askey--Wilson polynomials}
We will use the following standardization and
notation for \emph{Askey--Wilson polynomials}\,:
\begin{equation}
R_n[z]=R_n[z;a,b,c,d\,|\,q]:=
\qhyp43{q^{-n},q^{n-1}abcd,az,az^{-1}}{ab,ac,ad}{q,q}.
\label{30}
\end{equation}
These are symmetric Laurent polynomials of degree $n$ in $z$,
so they are ordinary polynomials of degree $n$ in $x:=\thalf(z+z^{-1})$.
The polynomials \eqref{30} are related to the Askey--Wilson polynomials
$p_n(x;a,b,c,d\,|\,q)$ in usual notation \cite[(1.15)]{9},
\cite[(14.1.1)]{2} by
\begin{equation}
R_n[z;a,b,c,d\,|\,q]
=\frac{a^n}{(ab,ac,ad;q)_n}\,
p_n\big(\thalf(z+z^{-1});a,b,c,d\,|\,q\big).
\label{31}
\end{equation}

If $|a|,|b|,|c|,|d|\le1$ such that pairwise products of $a, b, c, d$
are not equal
to 1 and such that non-real parameters occur in complex conjugate pairs,
then the Askey--Wilson polynomials are orthogonal with respect to a 
non-negative weight
function on $x=\thalf(z+z^{-1})\in[-1,1]$. For convenience we give
this orthogonality in the variable $z$ on the unit circle, where the integrand
is invariant under $z\to z^{-1}$:
\begin{equation}
\int_{|z|=1}R_m[z]\,R_n[z]\,
w[z]\,\frac{\dup z}{\iup z}
=h_n\,\de_{m,n},
\label{36}
\end{equation}
where
\begin{align}
w[z]=w[z;a,b,c,d;q]&
=\left|\frac{(z^2;q)_\iy}{(az,bz,cz,dz\,|\,q)_\iy}\right|^2,
\label{37}\sLP
h_0=h_0[a,b,c,d\,|\,q]&=\frac{4\pi (abcd;q)_\iy}{(q,ab,ac,ad,bc,bd,cd;q)_\iy},
\label{38}
\end{align}
and
where the explicit expression for $h_n$
can be obtained from
\cite[(14.1.2)]{2} together with \eqref{31}.
\subsection{Continuous \boldmath$q$-ultraspherical polynomials}
The continuous $q$-ultraspherical polynomials are a one-parameter
subfamily of the Askey--Wilson polynomials \eqref{30}.
For them we will use the following standardization and notation:
\begin{equation}
\begin{split}
R_n^{\be;q}[z]&=R_n^{\be;q}\big(\thalf(z+z^{-1})\big):=
R_n\big[z;q^{\frac14}\be^\half,q^{\frac34}\be^\half,
-q^{\frac14}\be^\half,-q^{\frac34}\be^\half\,|\,q\big]\\
&=\qhyp43{q^{-n},\be^2q^{n+1},q^{\frac14}\be^\half z,
q^{\frac14}\be^\half z^{-1}}{\be q,-\be q^\half,-\be q}{q,q}.
\end{split}
\label{10}
\end{equation}
The polynomials \eqref{10} are related to the
continuous $q$-ultraspherical polynomials in usual notation
\cite[\S14.10.1]{2} by
\begin{equation}
R_n^{\be;q}(x)=q^{\frac14 n}\be^{\half n}\,\frac{(q;q)_n}{(q\be^2;q)_n}\,
C_n\big(x;q^\half \be\,|\,q\big).
\label{96}
\end{equation}
The continuous $q$-ultraspherical polynomials with
$\be=q^\al$
tend to the ultraspherical polynomials \eqref{90} as $q\uparrow1$:
\begin{equation*}
\lim_{q\uparrow1} R_n^{q^\al;q}(x)=R_n^\al(x).
\end{equation*}

In view of  \cite[(3.10.13)]{3} we can represent $R_n^{\be;q}$
by a different $q$-hypergeometric expression:
\begin{equation}
R_n^{\be;q}[z]=
\qhyp43{q^{-\half n},q^{\half n+\half}\be,q^{\frac14}\be^\half z,
q^{\frac14}\be^\half z^{-1}}
{-q^\half\be,(q\be)^\half,-(q\be)^\half}{q^\half,q^\half}.
\label{42}
\end{equation}
In particular, 
\begin{equation*}
R_n^{\be;q}\big[q^{-\half m-\frac14}\be^{-\half}\big]
=\qhyp43{q^{-\half n},q^{\half n+\half}\be,q^{-\half m},q^{\half m+\half}\be}
{-q^\half\be,(q\be)^\half,-(q\be)^\half}{q^\half,q^\half}
\qquad(m,n=0,1,2,\ldots).
\end{equation*}
Hence we have the \emph{duality}
\begin{equation}
R_n^{\be;q}\big[q^{-\half m-\frac14}\be^{-\half}\big]
=R_m^{\be;q}\big[q^{-\half n-\frac14}\be^{-\half}\big]
\qquad(m,n=0,1,2,\ldots).
\label{44}
\end{equation}

Note the special value
\begin{equation*}
R_n^{\be;q}
\big[q^{\frac14}\be^\half\big]=1
\end{equation*}
and the coefficient of the term of highest degree
\begin{equation}
R_n^{\be;q}(x)=2^n (q^\half \beta)^{\half n}\,
\frac{(q^\half \be;q)_n}{(q\be^2;q)_n}\,x^n
+\mbox{terms of lower degree}.
\label{2}
\end{equation}

For $0<\be<q^{-\half}$
the polynomials $R_n^{\be;q}(x)$ are orthogonal on $[-1,1]$
with respect to the even weight function
\begin{equation}
w_{\be,q}(x):=(1-x^2)^{-\half} \left|\frac{(e^{2i\tha};q)_\iy}
{(q^\half\be e^{2i\tha};q)_\iy}\right|^2,\quad x=\cos\tha,
\label{1}
\end{equation}
see \cite[(14.10.18)]{2}.
This weight function satisfies the recurrence
\begin{align}
\frac{w_{q\be,q}(x)}{w_{\be,q}(x)}
&=\big(1+q^\half\be\big)^2-4q^\half\be x^2\nonu\\
&=4q^\half\be\big(a^2-x^2\big),\qquad
a=\thalf\big(q^{\frac14}\be^\half+q^{-\frac14}\be^{-\half}\big).\label{11}
\end{align}
We will need the difference formula
\begin{multline}
R_n^{\be;q}(x)-R_{n-2}^{\be;q}(x)
=\frac{4 q^{-\half n+\frac32}\be}{(1+q^\half\be)(1+q\be)}\,
\frac{1-q^{n-\half}\be}{1-q\be}\\
\times\Big(x^2
-\big(\thalf(q^{\frac14}\be^\half+q^{-\frac14}\be^{-\half})\big)^2\Big)R_{n-2}^{q\be;q}(x)\qquad
(n\ge2).\label{3}
\end{multline}
{\bf Proof of \eqref{3}.}\quad
More generally, let $w(x)=w(-x)$ be an even weight function on $[-1,1]$,
let $p_n(x)=k_nx^n+\cdots\;$ be orthogonal polynomials on $[-1,1]$ with
respect to the weight function $w(x)$, and let
$q_n(x)=k_n'x^n+\cdots\;$ be orthogonal polynomials on $[-1,1]$ with
respect to the weight function $w(x)(a^2-x^2)$ ($a\ge1$).
Assume that $p_n$ and $q_n$
are normalized by $p_n(a)=1=q_n(a)$. Let $n\ge2$.
Then $p_n(x)-p_{n-2}(x)$ vanishes for $x=\pm a$. Hence
$(p_n(x)-p_{n-2}(x))/(x^2-a^2)$ is a polynomial of degree $n-2$.
It is immediately  seen that $x^k$ ($k<n-2$) is orthogonal to this polynomial with respect to the weight function $w(x)(a^2-x^2)$ on $[-1,1]$.
We conclude that
\[
p_n(x)-p_{n-2}(x)=\frac{k_n}{k_{n-2}'}\,(x^2-a^2)q_{n-2}(x)\qquad(n\ge2).
\]
Now specialize to the weight function \eqref{1} and use
\eqref{11} and \eqref{2}.\qed
\bPP
From \cite[(14.10.17)]{2} we have
\begin{equation}
R_n^{\be;q}(\cos\tha)=q^{\frac14 n}\be^{\half n}\,
\frac{(q^\half \be;q)_n}{(q\be^2;q)_n}\,\eup^{\iup n\tha}
\qhyp21{q^{-n},q^\half \be}{q^{-n+\half} \be^{-1}}{q,q^\half\be^{-1}
\eup^{-2\iup n\tha}}.
\label{84}
\end{equation}
A limit case of \eqref{84} yields the \emph{continuous $q$-Hermite polynomials}
(see \cite[(14.26.1)]{2}):
\begin{equation}
H_n(\cos\tha\,|\,q):=\eup^{\iup n\tha}
\qhyp20{q^{-n},0}-{q,q^n \eup^{-2\iup \tha}}.
\label{86}
\end{equation}
So by \eqref{84} and \eqref{86} we have
\begin{equation}
H_n(x\,|\,q)=
q^{-\frac14 n}\,\lim_{\be\downarrow0}\,\be^{-\half n} R_n^{\be;q}(x).
\label{85}
\end{equation}
\subsection{\boldmath$q$-Racah polynomials}
We will consider \emph{$q$-Racah polynomials}
\cite[\S14.2]{2}
\begin{equation}
R_n(q^{-x}+\ga\de q^{x+1};\al,\be,\ga,\de\,|\,q):=
\qhyp43{q^{-n},q^{n+1}\al\be,q^{-x},q^{x+1}\ga\de}
{q\al,q\be\de,q\ga}{q,q}
\label{7}
\end{equation}
for $\ga=q^{-N-1}$, where $N\in\{1,2,\ldots\}$, and for
$n\in\{0,1,\ldots,N\}$.
They are discrete cases of the Askey--Wilson polynomials \eqref{30}.
The polynomials \eqref{7} are orthogonal polynomials on
the finite $q$-quadratic set
$\{q^{-x}+\ga\de q^{x+1}\mid x\in\{0,1,\ldots,N\}\}$:
\begin{equation}
\sum_{x=0}^N(R_mR_n)(q^{-x}+\ga\de q^{x+1};\al,\be,\ga,\de\,|\,q)\,
w_{\al,\be,\ga,\de;q}(x)=h_{n;\al,\be,\ga,\de;q}\,\de_{m,n}
\label{83}
\end{equation}
with
\begin{align}
w_{\al,\be,\ga,\de;q}(x)&:=\frac{1-\ga\de q^{2x+1}}{(\al\be q)^x (1-\ga\de q)}\,
\frac{(\al q,\be\de q,\ga q,\ga\de q;q)_x}
{(q,\al^{-1}\ga\de q,\be^{-1}\ga q,\de q;q)_x}\,,
\label{5}\sLP
\frac{h_{n;\al,\be,\ga,\de;q}}{h_{0;\al,\be,\ga,\de;q}}
&:=\frac{(1-\al\be q)(q\ga\de)^n}{1-\al\be q^{2n+1}}\,
\frac{(q,q\be,q\al\be\ga^{-1},q\al\de^{-1};q)_n}
{(q\al,q\al\be,q\ga,q\be\de;q)_n}\,,\label{20}
\\
h_{0;\al,\be,\ga,\de;q}&:=\sum_{x=0}^N w_{\al,\be,\ga,\de;q}(x)
=\frac{(q^2\al\be,\de^{-1};q)_N}{(q\al\de^{-1},q\be;q)_N}\qquad(\ga=q^{-N-1}).
\label{13}
\end{align}

Clearly $R_n(1+q^{-N}\de;\al,\be,q^{-N-1},\de\,|\,q)=1$ while,
by \eqref{7} and the $q$-Saalsch\"utz formula
\cite[(1.7.2)]{3}, we can evaluate
the $q$-Racah polynomial for $x=N$:
\begin{equation}
R_n(q^{-N}+\de;\al,\be,q^{-N-1},\de\,|\,q)=
\frac{(q\be,q\al\de^{-1};q)_n}{(q\al,q\be\de;q)_n}\,\de^n.
\label{8}
\end{equation}

The backward shift operator equation
\cite[(14.2.10)]{2} can be rewritten as
\begin{align}
&w_{\al,\be,\ga,\de;q}(x)
R_n\big(q^{-x}+\ga\de q^{x+1};\al,\be,\ga,\de\,|\,q\big)\nonu\\
&=\frac{1-q^2\ga\de}{q^{-x}-\ga\de q^{x+2}}\,w_{q\al,q\be,q\ga,\de;q}(x)\,
R_{n-1}\big(q^{-x}+\ga\de q^{x+2};q\al,q\be,q\ga,\de\,|\,q\big)\nonu\\
&-\frac{1-q^2\ga\de}{q^{-x+1}-\ga\de q^{x+1}}\,
w_{q\al,q\be,q\ga,\de;q}(x-1)\,
R_{n-1}\big(q^{-x+1}+\ga\de q^{x+1};q\al,q\be,q\ga,\de\,|\,q\big).
\label{4}
\end{align}
This holds for $x=1,\ldots,N$ while  for $x=0$ \eqref{4} remains true if
we put the second term on the right equal to 0. In the case $x=N$
the first term on the right is equal to zero because of \eqref{5},
and the identity \eqref{4} can be checked by using
\eqref{5} and \eqref{8}.

Hence, for a function $f$ on $\{0,1,\ldots,N\}$ we have
\begin{multline}
\sum_{x=0}^N
w_{\al,\be,\ga,\de;q}(x)
R_n\big(q^{-x}+\ga\de q^{x+1};\al,\be,\ga,\de\,|\,q\big)\,f(x)
=\sum_{x=0}^{N-1}
\frac{1-q^2\ga\de}{q^{-x}-\ga\de q^{x+2}}\\
\times w_{q\al,q\be,q\ga,\de;q}(x)\,
R_{n-1}\big(q^{-x}+\ga\de q^{x+2};q\al,q\be,q\ga,\de\,|\,q\big)\,
\big(f(x)-f(x+1)\big).
\label{9}
\end{multline}
\section{The dual addition formula for continuous
\boldmath$q$-ultraspherical\\ polynomials}
\label{75}
\subsection{The Rahman--Verma addition formula}
The $q$-analogue of the product formula for ultraspherical polynomials
\cite[(18.17.5)]{11} was given by Rahman \& Verma \cite[(1.20)]{10}.
It uses a different choice of parameter for the
$q$-ultraspherical polynomials:
\begin{equation}
R_n^{\textup{R-V}}[z]=R_n^{\textup{R-V};a;q}[z]:=R_n^{q^{-\half}a^2;q}[z],
\label{39}
\end{equation}
where we have our own notation \eqref{10} on the right.
Then the duality \eqref{44} takes the form
\begin{equation*}
R_n^{\textup{R-V}}\big[q^{-\half m} a^{-1}\big]=
R_m^{\textup{R-V}}\big[q^{-\half n} a^{-1}\big]
\qquad(m,n=0,1,2,\ldots)
\end{equation*}
or, in terms of special Askey--Wilson polynomials,
\begin{multline}
R_n\big[q^{-\half m} a^{-1};a,q^\half a,-a,-q^\half a\,|\,q\big]\\
=R_m\big[q^{-\half n} a^{-1};a,q^\half a,-a,-q^\half a\,|\,q\big]\qquad
(m,n=0,1,2,\ldots).
\label{47}
\end{multline}

In terms of the polynomials \eqref{39} and with usage of \eqref{37},
\eqref{38} the Rahman--Verma product formula reads as follows:
\begin{equation*}
R_n^{\textup{R-V}}[u]\,R_n^{\textup{R-V}}[v]=
\int_{|z|=1} R_n^{\textup{R-V}}[z]\,\frac{w[z;auv,au^{-1}v^{-1},auv^{-1},au^{-1}
v\,|\,q]}
{h_0(auv,au^{-1}v^{-1},auv^{-1},au^{-1}v\,|\,q)}\,\frac{\dup z}{\iup z}
\end{equation*}
with $|u|,|v|=1$, $0<a<1$.
This suggests an expansion
\begin{equation*}
R_n^{\textup{R-V}}[z]=\sum_{k=0}^n c_k\,
R_k\big[z;auv,au^{-1}v^{-1},auv^{-1},au^{-1}v\,|\,q\big],
\end{equation*}
where the term $c_0$ equals $R_n^{\textup{R-V}}[u]\,R_n^{\textup{R-V}}[v]$.
Indeed, \cite[(1.24)]{10} gives the addition formula
\begin{align}
R_n\big[z;a,q^\half a,-a&,-q^\half a\,|\,q\big]
=\sum_{k=0}^n\frac
{(-1)^k q^{\half k(k+1)}(q^{-n},a^2,q^n a^4,q^{-1}a^4;q)_k}
{(q,q^\half a^2,-q^\half a^2,-a^2;q)_k (q^{-1}a^4;q)_{2k}}\nonu\\
&\times u^{-k} (a^2u^2;q)_k\,
R_{n-k}\big[u;q^{\half k}a,q^{\half(k+1)}a,-q^{\half k}a,-q^{\half(k+1)}a
\,|\,q\big]\nonu\\
&\times v^{-k} (a^2v^2;q)_k\,
R_{n-k}\big[v;q^{\half k}a,q^{\half(k+1)}a,-q^{\half k}a,
-q^{\half(k+1)}a\,|\,q\big]\nonu\\
&\times R_k\big[z;auv,au^{-1}v^{-1},auv^{-1},au^{-1}v\,|\,q\big].
\label{45}
\end{align}
The addition formula \cite[(18.18.8)]{11} for ultraspherical polynomials
can be obtained as
limit case for $q\uparrow1$ of \eqref{45}.
\subsection{The dual addition formula}
As mentioned in \cite[(4.18)]{4},
Rogers already gave the linearization formula for continuous
$q$-ultraspherical polynomials in 1895.
Here we refer for this formula to \cite[(10.11.10)]{1}.
It can be written in notation \eqref{10} as
\begin{multline}
R_\ell^{\be;q}(x) R_m^{\be;q}(x)
=\frac{(q;q)_\ell(q;q)_m}{(q\be^2;q)_\ell(q\be^2;q)_m}
\sum_{j=0}^{\min(\ell,m)} \frac{1-q^{\ell+m-2j+\half}\beta}
{1-q^\half\beta}\;
\frac{(q^\half\beta;q)_j}{(q;q)_j}\\
\times\frac{(q^\half\beta;q)_{\ell-j}}{(q;q)_{\ell-j}}\,
\frac{(q^\half\beta;q)_{m-j}}{(q;q)_{m-j}}\,
\frac{(q\be^2;q)_{\ell+m-j}}{(q^{\frac32}\be;q)_{\ell+m-j}}\,
(q^\half\be)^j R_{\ell+m-2j}^{\be;q}(x).
\label{12}
\end{multline}
By the earlier assumption $0<\be<q^{-\half}$ the linearization
coefficients in \eqref{12} are nonnegative.

From now on assume without loss of generality that $\ell\ge m$.
Specialization of \eqref{5} and \eqref{13} gives
\begin{multline*}
w_{\be q^{-\half},\be q^{-\half},q^{-m-1},\be^{-1}q^{-\ell-\half};q}(j)
=\frac{(q^{\half}\be;q)_{\ell+m}}{(q\be^2;q)_{\ell+m}}\,
\frac{(q;q)_{\ell}}{(q^\half \be;q)_{\ell}}\,
\frac{(q;q)_{m}}{(q^\half \be;q)_{m}}\\
\times
\frac{1-q^{\ell+m-2j+\half}\beta}{1-q^{\half}\beta}\;
\frac{(q^\half\beta;q)_j}{(q;q)_j}\,
\frac{(q^\half \be;q)_{\ell-j}}{(q;q)_{\ell-j}}\,
\frac{(q^\half \be;q)_{m-j}}{(q;q)_{m-j}}\,
\frac{(q\be^2;q)_{\ell+m-j}}{(q^{\frac32}\be;q)_{\ell+m-j}}\,
(q^\half\be)^j
\end{multline*}
and
\begin{equation}
h_{0;\be q^{-\half},\be q^{-\half},q^{-m-1},\be^{-1}q^{-\ell-\half};q}
=\frac{(q\be^2;q)_\ell(q\be^2;q)_m}{(q\be^2;q)_{\ell+m}}\,
\frac{(q^\half\be;q)_{\ell+m}}{(q^\half\be;q)_\ell (q^\half\be;q)_m}\,.
\label{16}
\end{equation}
The linearization formula \eqref{12} can now be seen to have the
equivalent concise expression
\begin{equation}
R_\ell^{\be;q}(x) R_m^{\be;q}(x)=
\sum_{j=0}^m\frac
{w_{\be q^{-\half},\be q^{-\half},q^{-m-1},\be^{-1}q^{-\ell-\half};q}(j)}
{h_{0;\be q^{-\half},\be q^{-\half},q^{-m-1},\be^{-1}q^{-\ell-\half};q}}
\,R_{\ell+m-2j}^{\be;q}(x).
\label{14}
\end{equation}
This identity can be considered as giving the constant term
of an expansion of $R_{\ell+m-2j}^{\be;q}(x)$ as a function of $j$
in terms of $q$-Racah polynomials
\[
R_k(q^{-j}+\be^{-1} q^{j-\ell-m-\half};
\be q^{-\half},\be q^{-\half},q^{-m-1},\be^{-1}q^{-\ell-\half}\,|\,q).
\]
The general terms of this expansion will be obtained by evaluating
the sum
\begin{multline}
S_{k,\ell,m}^{\be;q}(x):=\sum_{j=0}^m
w_{\be q^{-\half},\be q^{-\half},q^{-m-1},\be^{-1}q^{-\ell-\half};q}(j)
\,R_{\ell+m-2j}^{\be;q}(x)\\
\times R_k(q^{-j}+\be^{-1} q^{j-\ell-m-\half};
\be q^{-\half},\be q^{-\half},q^{-m-1},\be^{-1}q^{-\ell-\half}\,|\,q),
\label{15}
\end{multline}
where we still assume $l\ge m$ and where $k\in\{0,\ldots,m\}$.
\begin{theorem}
The sum \eqref{15} can be evaluated as
\begin{multline}
S_{k,\ell,m}^{\be;q}(x)=
\frac{(q^{\half(\ell+m+1)}\be)^k
(\be^{-1} q^{-\ell-m+\half};q)_k}
{(-q^\half \be,\pm q\be;q)_k}\,
(\pm q^{\frac14}\be^\half z,\pm q^{\frac14}\be^\half z^{-1};q^\half)_k\\
\times \frac{(q^{2k+1}\be^2;q)_{\ell-k}(q^{2k+1}\be^2;q)_{m-k}}
{(q^{2k+1}\be^2;q)_{\ell+m-2k}}\,
\frac{(q^{k+\half}\be;q)_{\ell+m-2k}}
{(q^{k+\half}\be;q)_{\ell-k} (q^{k+\half}\be;q)_{m-k}}
R_{\ell-k}^{q^k\be;q}(x) R_{m-k}^{q^k\be;q}(x).
\label{18}
\end{multline}
Here we use the conventions that
$(\pm a;q)_n:=(a;q)_n (-a;q)_n$ and $x=\half(z+z^{-1})$.
\end{theorem}
\Proof
In \eqref{15} put
$f(j):=R_{\ell+m-2j}^{\be;q}(x)$.
Then comparison of \eqref{15} with \eqref{9} gives
\begin{align*}
&S_{k,\ell,m}^{\be;q}(x)=\sum_{j=0}^m
w_{\be q^{-\half},\be q^{-\half},q^{-m-1},\be^{-1}q^{-\ell-\half};q}(j)\\
&\qquad\qquad\qquad\times
R_k(q^{-j}+\be^{-1} q^{j-\ell-m-\half};
\be q^{-\half},\be q^{-\half},q^{-m-1},\be^{-1}q^{-\ell-\half}\,|\,q)\,f(j)\\
&\qquad\qquad=\sum_{j=0}^{m-1}
\frac{1-\be^{-1}q^{-\ell-m+\half}}{q^{-j}-\be^{-1}q^{-\ell-m+j+\half}}\,
w_{\be q^\half,\be q^\half,q^{-m},\be^{-1}q^{-\ell-\half};q}(j)\\
&\qquad\qquad\qquad\times R_{k-1}(q^{-j}+\be^{-1} q^{j-\ell-m+\half};
\be q^\half,\be q^\half,q^{-m},\be^{-1}q^{-\ell-\half}\,|\,q)
\big(f(j)-f(j+1)\big).
\end{align*}
We can handle the factor $f(j)-f(j+1)$ in the right part above by
using \eqref{3}:
\begin{align*}
&f(j)-f(j+1)
=R_{\ell+m-2j}^{\be;q}(x)-R_{\ell+m-2j-2}^{\be;q}(x)\\
&\quad=\frac{4\be^2 q^{\half\ell+\half m+1}
(q^{-j}-\be^{-1}q^{-\ell-m+j+\half})}
{(1+q^\half\be)(1-q^2\be^2)}\,
\Big(\tfrac14
\big(q^{\frac14}\be^\half+q^{-\frac14}\be^{-\half}\big)^2-x^2\Big)\,
R_{\ell+m-2j-2}^{q\be;q}(x).
\end{align*}
So, with $x=\thalf(z+z^{-1})$,
\begin{align*}
&S_{k,\ell,m}^{\be;q}(x)
=\frac{4\be q^{-\half\ell-\half m+\frac32}
(1-\be q^{\ell+m-\half})}
{(1+q^\half\be)(1-q^2\be^2)}\\
&\qquad\times\Big(\tfrac14
\big(q^{\frac14}\be^\half+q^{-\frac14}\be^{-\half}\big)^2-x^2\Big)
\sum_{j=0}^{m-1}
w_{\be q^\half,\be q^\half,q^{-m},\be^{-1}q^{-\ell-\half};q}(j)\\
&\qquad\qquad\times R_{k-1}(q^{-j}+\be^{-1} q^{j-\ell-m+\half};
\be q^\half,\be q^\half,q^{-m},\be^{-1}q^{-\ell-\half}\,|\,q)
R_{\ell+m-2j-2}^{q\be;q}(x)\\
&\quad=\frac{q^{\half\ell+\half m+\half}\be(1-\be^{-1} q^{-l-m+\half})}
{(1+q^\half\be)(1-q^2\be^2)}\,
(1+q^{\frac14}\be^\half z)(1-q^{\frac14}\be^\half z)
(1+q^{\frac14}\be^\half z^{-1})(1-q^{\frac14}\be^\half z^{-1})\\
&\qquad\qquad\qquad\qquad\qquad\qquad\qquad\qquad
\times S_{k-1,\ell-1,m-1}^{q\be,q}(x).
\end{align*}
Iteration gives
\begin{equation}
S_{k,\ell,m}^{\be;q}(x)
=
\frac{(q^{\half(\ell+m+1)}\be)^k
(\be^{-1} q^{-\ell-m+\half};q)_k}
{(-q^\half \be,\pm q\be;q)_k}\,
(\pm q^{\frac14}\be^\half z,\pm q^{\frac14}\be^\half z^{-1};q^\half)_k\,
S_{0,\ell-k,m-k}^{q^k\be;q}(x).
\label{17}
\end{equation}
By \eqref{15}
\begin{equation}
S_{0,\ell,m}^{\be;q}(x)=
h_{0;\be q^{-\half},\be q^{-\half},q^{-m-1},\be^{-1}q^{-\ell-\half};q}
R_\ell^{\be;q}(x) R_m^{\be;q}(x).
\label{19}
\end{equation}
Hence, by \eqref{16},
\begin{multline*}
S_{0,\ell-k,m-k}^{q^k\be;q}(x)=
h_{0;\be q^{k-\half},\be q^{k-\half},q^{k-m-1},\be^{-1}q^{-\ell-\half};q}
R_{\ell-k}^{q^k\be;q}(x) R_{m-k}^{q^k\be;q}(x)\\
=\frac{(q^{2k+1}\be^2;q)_{\ell-k}(q^{2k+1}\be^2;q)_{m-k}}
{(q^{2k+1}\be^2;q)_{\ell+m-2k}}\,
\frac{(q^{k+\half}\be;q)_{\ell+m-2k}}
{(q^{k+\half}\be;q)_{\ell-k} (q^{k+\half}\be;q)_{m-k}}
R_{\ell-k}^{q^k\be;q}(x) R_{m-k}^{q^k\be;q}(x).
\end{multline*}
Substitution of this last result in \eqref{17} yields \eqref{18}.\qed
\begin{remark}
\label{99}
An integrated form of \eqref{18} is the same as a special case
of the formula given in \cite[Remark 6.5]{46} (corrected version of
\cite[Remark 6.5]{39}).
Indeed, in \eqref{18} rewrite the three continuous $q$-ultraspherical
polynomials in the standard notation \eqref{96}, replace $\ell$ by
$n$ and $\be$ by $q^{-\half}\be$, multiply both sides by
$C_{m+n-2t}(x;\be\,|\,q)$ ($0\le t\le m$)
times its weight function, and integrate both sides
over $x\in [-1,1]$ (see \cite[(14.10.18)]{2}). Then write the Racah
polynomial on the \RHS\
by \eqref{7} as a balanced ${}_4\phi3$ and apply to this
Sears' transformation \cite[(III.15)]{3}
(with $n$, $a$, $b$, $c$, $d$, $e$, $f$ replaced by
$k$, $q^{k-1}\be^2$, $q^{-t}$, $q^{-n-m+t}\be^{-1}$, $q^{-m}$, $q^{-n}$, $\be$).
We arrive at the formula in \cite[Remark 6.5]{46} for $\al=q^{-1}\be$.

Note that \eqref{18} is not equivalent to its integrated forms if we consider
these only for $0\le t\le m$. The \RHS\ of \eqref{18} is a polynomial of
degree $\ell+m$ in $x$, so we have to consider also the integrals for
$m<t\le\half(\ell+m)$, of which we know a priori that they vanish. But
\cite[Remark 6.5]{46} does not consider integrals for $t>m$. However,s in
\cite[Lemma 6.4]{46} (the case $\be=1$, $\al=0$ of \cite[Remark 6.5]{46})
the integral for $t>m$ is stated to be zero.

As observed in \cite[end of \S4]{8}, integrated forms of the $q=1$ limit
\cite[(4.5)]{8} of \eqref{18} coincide with special cases of
\cite[(2.6)]{40}.
\end{remark}
\begin{theorem}[Dual addition formula]
For $j\in\{0,\ldots,m\}$ there is the expansion
\begin{align}
R_{\ell+m-2j}^{\be;q}(x)&=\sum_{k=0}^{\min(l,m)} q^{\half k(k+\ell+m+2)}\be^k\,
\frac{1-\be^2 q^{2k}}{1-\be^2 q^k}\,
\frac{(q^{-\ell},q^{-m},q\be^2;q)_k}{(q\be,q\be,q;q)_k}\nonu\\
&\qquad\times\frac{\prod_{i=0}^{k-1}
\big(4q^{i+\half}\be x^2-(1+q^{i+\half}\be)^2\big)}
{(-q^\half \be;q^\half)_{2k}^2}\,
R_{\ell-k}^{q^k\be;q}(x) R_{m-k}^{q^k\be;q}(x)\nonu\\
&\qquad\times R_k(q^{-j}+\be^{-1} q^{j-\ell-m-\half};
\be q^{-\half},\be q^{-\half},q^{-m-1},\be^{-1}q^{-\ell-\half}\,|\,q).
\label{46}
\end{align}
\end{theorem}
\Proof
Assume $l\ge m$.
By \eqref{17} and \eqref{19}
\begin{multline*}
S_{k,\ell,m}^{\be;q}(x)
=\frac{(-1)^k q^{\half k(k-\ell-m+1)}}
{(-q^\half\be;q)_k^2(q^2\be^2;q^2)_k^2}\,
\frac{(q\be^2;q)_{\ell+k} (q\be^2;q)_{m+k} (q^\half\be;q)_{\ell+m}}
{(q^\half\be;q)_{\ell+m}(q\be^2;q)_\ell(q\be^2;q)_m}\\
\times h_{0;\be q^{-\half},\be q^{-\half},q^{-m-1},\be^{-1}q^{-\ell-\half};q}
(\pm q^{\frac14}\be^\half z,\pm q^{\frac14}\be^\half z^{-1};q^\half)_k
R_{\ell-k}^{q^k\be,q}(x) R_{m-k}^{q^k\be,q}(x).
\end{multline*}
By Fourier-$q$-Racah inversion we obtain
\begin{align*}
&R_{\ell+m-2j}^{\be,q}(x)=\sum_{k=0}^m
\frac{(-1)^k q^{\half k(k-\ell-m+1)}}
{(-q^\half\be;q)_k^2(q^2\be^2;q^2)_k^2}\,
\frac{(q\be^2;q)_{\ell+k} (q\be^2;q)_{m+k} (q^\half\be;q)_{\ell+m}}
{(q^\half\be;q)_{\ell+m}(q\be^2;q)_\ell(q\be^2;q)_m}\\
&\qquad\times\frac
{h_{0;\be q^{-\half},\be q^{-\half},q^{-m-1},\be^{-1}q^{-\ell-\half};q}}
{h_{k;\be q^{-\half},\be q^{-\half},q^{-m-1},\be^{-1}q^{-\ell-\half};q}}\,
(\pm q^{\frac14}\be^\half z,\pm q^{\frac14}\be^\half z^{-1};q^\half)_k
R_{\ell-k}^{q^k\be;q}(x) R_{m-k}^{q^k\be;q}(x)\\
&\qquad\times
R_k(q^{-j}+\be^{-1} q^{j-\ell-m-\half};
\be q^{-\half},\be q^{-\half},q^{-m-1},\be^{-1}q^{-\ell-\half}\,|\,q).
\end{align*}
Now use \eqref{20}.\qed
\bPP
If we put $\be=q^\al$ in  \eqref{46} and take the limit for $q\uparrow1$
then we arrive at the dual addition formula \eqref{29} for
ultraspherical polynomials.

For $j=0$ \eqref{46} takes the form
\begin{align}
R_{\ell+m}^{\be;q}(x)&=\sum_{k=0}^{\min(l,m)} q^{\half k(k+\ell+m+2)}\be^k\,
\frac{1-\be^2 q^{2k}}{1-\be^2 q^k}\,
\frac{(q^{-\ell},q^{-m},q\be^2;q)_k}{(q\be,q\be,q;q)_k}\nonu\\
&\qquad\times\frac{\prod_{i=0}^{k-1}
\big(4q^{i+\half}\be x^2-(1+q^{i+\half}\be)^2\big)}
{(-q^\half \be;q^\half)_{2k}^2}\,
R_{\ell-k}^{q^k\be;q}(x) R_{m-k}^{q^k\be;q}(x).
\label{98}
\end{align}
It has a similar structure as \cite[Exercise 8.12]{3}.
However, the formula there expands $R_{\ell+m}^{\be;q}(x)$ in terms of
$R_{\ell-k}^{\be;q}(x)R_{m-k}^{\be;q}(x)$ ($k = 0, 1, \ldots, \min(\ell, m)$).
A variant of \eqref{87} given by the specialization
of \cite[(9.4)]{42} to the case of continuous $q$-ultraspherical polynomials
is also essentially different from our formula.

Just as with  \eqref{29},
formula \eqref{98} can be rewritten as a matrix decomposition $S=LDU$ with $S$
symmetric, $L$ lower triangular, its transpose $U=L^\intercal$ upper
triangular and $D$ diagonal.
Aldenhoven \cite[Theorem 1.1]{37} earlier gave the inverse of the matrix $L$.
\subsection{A second proof of the dual addition formula}
We will now show that the addition formula \eqref{45} and the
dual additon formula \eqref{46} coincide when both formulas are
suitably restricted in their $x$ or $z$ variable.
This will follow from the duality \eqref{47}.

In \eqref{46} put $\be=q^{-\half} a^2$, $x=\thalf(z+z^{-1})$
and use \eqref{39}.
Then the dual addition formula takes the form
\begin{align}
&R_{\ell+m-2j}\big[z;a,q^\half a,-a,-q^\half a\,|\,q\big]\nonu\\
&=\sum_{k=0}^m (-1)^k q^{\half k(k+\ell+m+1)}a^{2k}\,
\frac{1-a^4 q^{2k-1}}{1-a^4 q^{k-1}}\,
\frac{(q^{-\ell},q^{-m},a^4;q)_k}{(q^\half a^2,q^\half a^2,q;q)_k}\,
\frac{(a^2z^2,a^2z^{-2};q)_k}
{(-a^2;q^\half)_{2k}^2}\nonu\\
&\quad\times
R_{\ell-k}\big[z;q^{\half k}a,q^{\half(k+1)}a,-q^{\half k}a,
-q^{\half(k+1)}a\,|\,q\big]\,
R_{m-k}\big[z;q^{\half k}a,q^{\half(k+1)}a,-q^{\half k}a,
-q^{\half(k+1)}a\,|\,q\big]\nonu\\
&\quad\times R_k\big(q^{-j}+q^{j-\ell-m}a^{-2};
q^{-1}a^2,q^{-1}a^2,q^{-m-1},q^{-\ell}a^{-2}\,|\,q\big).
\label{48}
\end{align}
Since both sides of \eqref{48} are symmetric Laurent polynomials in $z$,
verification of the identity for $z=q^{-\half n}a^{-1}$
($n=m,m+1,m+2,\ldots$)
will settle the identity for all $z$. Thus put $z=q^{-\half n}a^{-1}$
in \eqref{48} and use the duality \eqref{47} in the
polynomials $R_{\ell+m-2j}$, $R_{\ell-k}$ and $R_{m-k}$ occurring in
\eqref{48}. Furthermore, use \eqref{7} and \eqref{30}
in order to substitute
\begin{align*}
&R_k(q^{-j}+q^{j-\ell-m}a^{-2};
q^{-1}a^2,q^{-1}a^2,q^{-m-1},q^{-\ell}a^{-2}\,|\,q)\\
&\qquad=\qhyp43{q^{-k},q^{k-1}a^k,q^{-j},q^{j-\ell-m}a^{-2}}
{a^2,q^{-\ell},q^{-m}}{q,q}\\
&\qquad=R_k\big[q^{-\half(\ell+m-2j)}a^{-1};q^{-\half(\ell+m)}a^{-1},
q^{\half(\ell+m)}a^3,q^{\half(\ell-m)}a,q^{\half(m-\ell)}a\,|\,q\big].
\end{align*}
We obtain
\begin{align}
&R_n\big[q^{-\half(\ell+m-2j)}a^{-1};a,q^\half a,-a,-q^\half a\,|\,q\big]\nonu\\
&\qquad=\sum_{k=0}^n (-1)^k q^{\half k(k+\ell+m+1)}a^{2k}\,
\frac{1-a^4 q^{2k-1}}{1-a^4 q^{k-1}}\,
\frac{(q^{-\ell},q^{-m},a^4;q)_k}{(q^\half a^2,q^\half a^2,q;q)_k}\,
\frac{(q^{-n},q^n a^4;q)_k}{(-a^2;q^\half)_{2k}^2}\nonu\\
&\qquad\qquad\times
R_{n-k}\big[q^{-\half\ell}a^{-1};q^{\half k}a,q^{\half(k+1)}a,
-q^{\half k}a,-q^{\half(k+1)}a\,|\,q\big]\nonu\\
&\qquad\qquad\times
R_{n-k}\big[q^{-\half m}a^{-1};q^{\half k}a,q^{\half(k+1)}a,
-q^{\half k}a,-q^{\half(k+1)}a\,|\,q\big]\nonu\\
&\qquad\qquad\times
R_k\big[q^{-\half(\ell+m-2j)}a^{-1};q^{-\half(\ell+m)}a^{-1},
q^{\half(\ell+m)}a^3,q^{\half(\ell-m)}a,q^{\half(m-\ell)}a\,|\,q\big].
\label{49}
\end{align}
Because of the factor $(q^{-m};q)_k$ on the \RHS\ and since $n\ge m$
there was no harm to replace $m$ by $n$ as the upper bound of the
summation.

On the other hand, for integers $j,m,\ell$ such that
$0\le j\le m\le\ell$ and $m\le n$,
substitute $z=q^{-\half(\ell+m-2j)}a^{-1}$,
$u=q^{-\half\ell} a^{-1}$, $v=q^{-\half m}a^{-1}$ in \eqref{45} in order
to obtain
\begin{align}
&R_n\big[q^{-\half(\ell+m-2j)}a^{-1};a,q^\half a,-a,
-q^\half a\,|\,q\big]\nonu\\
&\qquad=\sum_{k=0}^n\frac
{(-1)^k q^{\half k(k+\ell+m+1)}a^{2k}(q^{-n},q^{-\ell},q^{-m},a^2,q^n a^4,q^{-1}a^4;q)_k}
{(q,q^\half a^2,-q^\half a^2,-a^2;q)_k (q^{-1}a^4;q)_{2k}}\nonu\\
&\qquad\qquad\times
R_{n-k}\big[q^{-\half\ell} a^{-1};q^{\half k}a,q^{\half(k+1)}a,-q^{\half k}a,-q^{\half(k+1)}a
\,|\,q\big]\nonu\\
&\qquad\qquad\times
R_{n-k}\big[q^{-\half m}a^{-1};q^{\half k}a,q^{\half(k+1)}a,-q^{\half k}a,-q^{\half(k+1)}a
\,|\,q\big]\nonu\\
&\qquad\qquad\times
R_k\big[q^{-\half(\ell+m-2j)}a^{-1};q^{-\half(\ell+m)}a^{-1},
q^{\half(\ell+m)}a^3,q^{\half(\ell-m)}a,q^{\half(m-\ell)}a\,|\,q\big].
\label{50}
\end{align}
An easy computation shows that \eqref{49} can be rewritten as
\eqref{50}. Thus we have shown that the addition formula \eqref{45}
implies the dual addition formula \eqref{46}. 
\section{A limit to continuous $q$-Hermite polynomials}
\label{62}
This section gives the $q$-analogue of the results in \cite[\S5]{8}. The
treatment given here is completely parallel to the one given there.

We will do a rescaling in the dual addition formula \eqref{46} such
that we can take the limit for $\be\downarrow0$. For this purpose observe that
the $q$-Racah polynomial \eqref{7}
has limits
\begin{equation}
\begin{split}
\lim_{\be\downarrow0}\be^j R_n\big(q^{-j}+\be^{-1}q^{-m-l+j-\half};\,
\be q^{-\half},\be q^{-\half},q^{-m-1},\be^{-1}q^{-l-\half};q\big)
&=\frac{(q^{-n};q)_j}{(q^{-l},q^{-m};q)_j}\,q^{(j-m-l-\half)j},\\
\lim_{\be\downarrow0}\be^n R_n\big(q^{-j}+\be^{-1}q^{-m-l+j-\half};\,
\be q^{-\half},\be q^{-\half},q^{-m-1},\be^{-1}q^{-l-\half};q\big)
&=\frac{(q^{-j};q)_n}{(q^{-l},q^{-m};q)_n}\,q^{(j-m-l-\half)n},
\end{split}
\label{79}
\end{equation}
where $l,m\ge \max(j,n)$.
Otherwise said,
\[
R_n\big(q^{-j}+\be^{-1}q^{-m-l+j-\half};\,\be q^{-\half},\be q^{-\half},
q^{-m-1},\be^{-1}q^{-l-\half};q\big)=
O(\be^{-\min(n,j)})
\]
as $\be\downarrow0$ with the order constant given in \eqref{79}.

Now, in \eqref{46}, multiply both sides
by $\be^{-\half(l+m-2j)}$ and let $\be\downarrow0$.
By \eqref{85} and \eqref{79} we obtain for $l\ge m$ that
\begin{multline}
q^{(l+m-j)j} (q^{-l},q^{-m};q)_j H_{l+m-2j}(x\,|\,q)=
\sum_{k=j}^m (-1)^k q^{k(l+m)} q^{-\half k(k-1)} (q^{-l},q^{-m};q)_k\\
\times H_{l-k}(x\,|\,q)H_{m-k}(x\,|\,q)\,
\frac{(q^{-k};q)_j}{(q;q)_k}\,,
\label{89}
\end{multline}
which may be called the dual addition formula for continuous
$q$-Hermite polynomials.
When writen equivalently as
\begin{align*}
H_{l+m-2j}(x\,|\,q) = \sum_{k=j}^m (-1)^{k-j} q^{(k-j)(l+m-2j+\half)}
q^{-\half(k-j)^2}&\\
\times\frac{(q^{-l+j},q^{-m+j};q)_{k-j}}{(q;q)_{k-j}}\,
&H_{l-k}(x\,|\,q)H_{m-k}(x\,|\,q),
\end{align*}
it is seen to be equivalent to its special case $j=0$, which can be written as
\begin{equation}
H_{l+m}(x\,|\,q) = \sum_{k=0}^{\min(l,m)} (-1)^{k}
q^{\half k(k-1)} (q;q)_k \qbinom lkq \qbinom mkq H_{l-k}(x\,|\,q)
H_{m-k}(x\,|\,q),
\label{87}
\end{equation}
where $\dstyle\qbinom lkq:=\frac{(q;q)_l}{(q;q)_k (q;q)_{l-k}}$ is a
$q$-binomial coefficient.
Formula \eqref{87} was first given, with two different proofs, by
Carlitz \cite[(1.8)]{43}, \cite[(3)]{44}. Here our $H_n(x\,|\,q)$ is related
to Carlitz' $H_n(z)$ by
$H_n(\half(z+z^{-1})\,|\,q)=z^{-n}H_n(z^2)$.
An (essentially different) variant of \eqref{87} is given in
\cite[(9.5)]{42}.

The $q=1$ limit of \eqref{87} is \cite[10.13(36)]{35},
which goes back to Nielsen (1918) and Burchnall (1941), see \cite[\S1]{42} 
for historical details.

Note that \eqref{87} gives a matrix decomposition $S=LDU$ with $S$
symmetric, $L$ lower triangular, its transpose $U=L^\intercal$ upper
triangular and $D$ diagonal. 
Aldenhoven \cite[Corollary 5.2]{37} earlier gave the inverse of the matrix $L$.

Next we want to consider the limit as $\be\downarrow0$ of \eqref{18} with
$S_{k,\ell,m}^{\be;q}$ given by \eqref{15}. Recall that \eqref{18} together with
\eqref{15} is
the dual of \eqref{46} in the sense of Fourier-$q$-Racah inversion.
Observe from \eqref{5}--\eqref{13} that
\begin{align}
\lim_{\be\downarrow0}\,\be^{-j}\,
w_{\be q^{-\half},\be q^{-\half},q^{-m-1},\be^{-1}q^{-l-\half};q}(j)
&=\frac{(q^{-l},q^{-m};q)_j}{(q;q)_j}\,q^{(l+m-j+\frac32)j},\label{80}\\
\lim_{\be\downarrow0}\,\be^n\,
h_{n;\be q^{-\half},\be q^{-\half},q^{-m-1},\be^{-1}q^{-l-\half};q}
&=\frac{(q;q)_n}{(q^{-l},q^{-m};q)_n}\,q^{-(l+m+\half)n}.\label{81}
\end{align}
In \eqref{18} multiply both sides by $\be^{-\half(l+m)+k}$
and let $\be\downarrow0$.
By \eqref{85}, \eqref{79} and \eqref{80}
we obtain for $l\ge m$ that
\begin{multline}
\sum_{j=k}^m q^{j(l+m)} q^{-j^2} (q^{-l},q^{-m};q)_j H_{l+m-2j}(x\,|\,q)\,
\frac{(q^{-j};q)_k q^{j(k+1)}}{(q;q)_j}\\
=(-1)^k q^{k(l+m)} q^{-\half k(k-1)} (q^{-l},q^{-m};q)_k
H_{l-k}(x\,|\,q) H_{m-k}(x\,|\,q).
\label{88}
\end{multline}
When written equivalently as
\begin{align*}
\sum_{j=k}^m q^{(j-k)(l+m-2k)} q^{-(j-k)(j-k-1)}\,
\frac{(q^{-(l-k)},q^{-(m-k)};q)_{j-k}}{(q;q)_{j-k}}\,H_{l+m-2j}(x\,|\,q)\\
=H_{l-k}(x\,|\,q) H_{m-k}(x\,|\,q)
\end{align*}
it can be seen, just as with \eqref{89}, to be equivalent 
to its special case $k=0$, which can be written as
\begin{equation}
\sum_{j=0}^{\min(l,m)} (q;q)_j \qbinom ljq \qbinom mjq H_{l+m-2j}(x\,|\,q)
=H_l(x\,|\,q) H_m(x\,|\,q).
\end{equation}
This is the linearization formula for continuous $q$-Hermite polynomials,
see \cite[(10.11.17)]{1}.

Just as with \eqref{46} and \eqref{18}, the identities  \eqref{89} and
\eqref{88} can be obtained from each other by a Fourier type inversion.
This no longer involves an orthogonal system as the $q$-Racah polynomials but
a biorthogonal system implied by
the \emph{biorthogonality relation}
\begin{equation}
\sum_{j=0}^\iy
\frac{(q^{-n};q)_j}{(q;q)_j}\,\frac{(q^{-j};q)_k\,q^{j(k+1)}}{(q;q)_k}
=\de_{k,n}
\label{82}
\end{equation}
(see Carlitz \cite[Theorem 2]{33} or
Krattenthaler \cite[(1.2)]{34} for $a_j=1$, $b_j=0$).
Note that the above sum in fact runs from $j=k$ to $n$.
For $k<n$ formula \eqref{82} is also equivalent to
${}_1\phi_0(q^{k-n};-;q,q)=
\sum_{j=0}^{n-k}\frac{(q^{k-n};q)_j}{(q;q)_j}\,q^j=0$.
\\[\smallskipamount]\indent
The biorthogonality \eqref{82} is also a limit case of
the $q$-Racah orthogonality relation \eqref{83}. Indeed, replace
$\al,\be,\ga,\de$ by
$\be q^{-\half},\be q^{-\half},q^{-m-1},\allowbreak\be^{-1}q^{-l-\half}$,
multiply both sides of \eqref{83} by $\be^n$, let $\be\downarrow0$,
and use \eqref{79}, \eqref{80} and \eqref{81}.
\paragraph{Acknowledgement}
The results of Section \ref{75} appeared earlier in
\emph{Proceedings of the 17th Annual Conference SSFA} (India) \cite{37},
a volume which is not widely distributed outside India
and is not online available.
That contribution was based on the 2017 R.~P. Agarwal Memorial Lecture, which
was delivered by the author at the 17th Annual
Conference SSFA in Bikaner, Rajasthan, India.
I thank Prof.\ M.~A. Pathan and Prof.\ S.~A. Ali for the invitation
to deliver this lecture. I thank Prof.\ M.~A. Pathan for permission to
republish part of this paper in afapted form here.

I am grateful to Erik Koelink for calling my attention to the formula
in \cite[Remark 6.5]{39} (which led to Remark \ref{99})
and to the discussion of Nielsen--Burchnall type
formulas in \cite{42}.

\begin{small}
\begin{quote}
T. H. Koornwinder, Korteweg-de Vries Institute, University of
Amsterdam,\\
P.O.\ Box 94248, 1090 GE Amsterdam, The Netherlands;\\
email: \texttt{thkmath@xs4all.nl}
\end{quote}
\end{small}
%

\end{document}